\documentclass[12pt]{article}
\usepackage[latin1]{inputenc}
\usepackage{amssymb}
\usepackage{amsmath}
\usepackage{amsfonts}
\usepackage{graphicx}
\usepackage{graphics}
\newtheorem{theorem}{Theorem}

\newtheorem{remark}[theorem]{Remark}
\newcommand{\R}{\mathbb{R}}
\newcommand{\po}{{\hspace*{-1ex}}{\bf .  }}

\def\qed{\ifhmode\unskip\nobreak\fi\ifmmode\ifinner\else
\hskip5 pt \fi\fi\hbox{\hskip5 pt \vrule width4 pt
height6 pt depth1.5 pt \hskip 1pt }}
\begin{document}

\title{Constant mean curvature hypersurfaces with single valued projections on planar domains}
\author{M. Dajczer and J. Ripoll}
\date{}
\maketitle

\begin{abstract}
A classical problem in constant mean curvature hypersurface theory is, for given $H\geq 0$,
to determine whether a compact submanifold $\Gamma^{n-1}$ of codimension two  
 in Euclidean space $\R_+^{n+1}$, having a single valued orthogonal
projection on $\R^n$, is the boundary of a graph with constant mean
curvature $H$ over a domain in $\R^n$. A well known result of Serrin
gives a sufficient condition, namely, $\Gamma$ is contained in a right
cylinder $C$ orthogonal to $\R^n$ with inner mean curvature
$H_C\geq H$. In this paper, we prove existence and uniqueness if the
orthogonal projection $L^{n-1}$ of $\Gamma$ on $\R^n$ has mean
curvature $H_L\geq-H$ and $\Gamma$ is contained in a cone $K$ with basis in
$\R^n$ enclosing a domain in $\R^n$ containing $L$ such
that the mean curvature of $K$ satisfies $H_K\geq H$. Our condition reduces
to Serrin's when the vertex of the cone is infinite.
\end{abstract}

\bigskip

A classical problem in the theory of constant mean curvature (CMC)
hypersurfaces is to determine whether an $(n\!-\!1)$-dimensional compact
embedded connected submanifold $\Gamma^{n-1}$ of Euclidean space
$\R_{+}^{n+1}=\{x_{n+1}\geq0\}$ with a single valued orthogonal
projection $\gamma^{n-1}$ on the hyperplane $\R^n=\{x_{n+1}=0\}$ is
the boundary of an $n$-dimensional graph  with given constant
mean curvature $H\ge0$ (called $H\!$-graph) over the domain enclosed by 
the submanifold $\gamma$. 

In the minimal case, it was shown by Finn \cite{F1}, \cite{F2} for surfaces
and then by Jenkins-Serrin \cite{JS} and Bakel'man \cite{B1}, \cite{B2} for
any dimension that there exists a unique $H\!$-graph if the projection of
$\Gamma$ onto $\R^n$ bounds a convex domain. For $H>0$, it was proved by  
Serrin \cite{Se} that a unique $H\!$-graph exists if the mean curvature $H_C$ of
the right cylinder $C(\Gamma)$ over $\Gamma$ orthogonal to $\R^n$
satisfies $H_C\geq H$ (considering the non-normalized mean curvature taken
with respect to the normal direction pointing to the simply connected
component of $\R^{n+1}\backslash C(\Gamma))$.

We observe that the $H$-convexity assumption of the domain in the above 
results can not be dispensed. For instance, it follows from the maximum principle
that a circle in $\R^2$ with radius strictly larger than $1/H$ cannot be the boundary
of an $H\!$-graph.
Nevertheless, the $H$-convexity hypothesis may be weakened. 
In fact, here we prove that Serrin's assumption is a special case of a more
general condition that applies to a large class of domains and boundary data.

To state our results, we first introduce some terminology.
Let $V\in\R_{+}^{n+1}$ be a point and $\gamma^{n-1}$ a compact smooth embedded
submanifold of $\R^n$. We denote by $D_{\gamma}$ the closure of the
domain enclosed by $\gamma$ and by $K_V(\gamma)$ the cone with base $\gamma$
and vertex $V$, i.e.,
$$
K_V(\gamma)=\{tV+(1-t)p:p\in\gamma\mbox{ and }t\in[0,1]\}.
$$
We refer to the right vertical cylinder 
$$
C(\gamma)=\{p+te_{n+1}:p\in\gamma\mbox{ and }t\in[0,+\infty)\}
$$ 
as the cone $K_{\infty}(\gamma)$ over
$\gamma$ with vertex at infinity. Then, we say that $K_V(\gamma)$ is an
$H\!$-cone for $V\in\R^{n+1}\cup\{\infty\}$ if either 
$K_V(\gamma)\backslash\{V\}$ or $K_{\infty}(\gamma)$ is a smooth
hypersurface having (non-normalized) mean curvature $H_K\geq H>0$ with
respect to the inner orientation.

Serrin's result may now be restated as follows: If $K_{\infty}(\gamma)$ is an
$H\!$-cone and $\Gamma$ a compact embedded hypersurface of $K_{\infty}
(\gamma)$ having a single projection onto $\gamma$, then there is a unique
$H\!$-graph with boundary $\Gamma$ over the domain enclosed by $\gamma$. Here,
we show that if $\Gamma\subset K_V(\gamma)$ is a compact embedded
hypersurface in an $H\!$-cone having a single orthogonal projection onto a
hypersurface $L\subset D_{\gamma}$ with mean curvature $H_L\ge-H$, then
there is a constant mean curvature $H\!$-graph with boundary $\Gamma$ over the
domain enclosed by $L$.
\medskip

We treat separately the cases where $\Gamma$ is smooth ($\mathcal{C}
^{2,\alpha})$ or only continuous and obtain results unifying the two
statements above in both cases. We first state the result for the smooth case.

\begin{theorem}\po\label{hcons smooth} 
Let $\Gamma^{n-1}\subset\R_{+}^{n+1}$ be a compact embedded connected 
$\mathcal{C}^{2,\alpha}$ 
submanifold having a single orthogonal projection onto a hypersurface
$L^{n-1}\subset\R^n$ such that $H_L\geq-H$ for some constant
$H>0$. If there is an $H\!$-cone $K_V(\gamma)\subset\R^{n+1}$ such
that $\Gamma\subset K_V(\gamma)$ and $L\subset D_{\gamma}$ if $V$ is finite,
then there is a unique $H\!$-graph of class $\mathcal{C}
^{2,\alpha}$  with boundary
$\Gamma$ over the domain enclosed by $L$.
\end{theorem}

In the case of continuous boundary we have the following result.

\begin{theorem}\po\label{hcons c0} 
Let $\Gamma^{n-1}\subset \R_{+}^{n+1}$ be a compact embedded connected 
$\mathcal{C}^0$ submanifold
having a single orthogonal projection onto a $C^2$ hypersurface
$L^{n-1}\subset\R^n$ such that $H_L\geq-H$ for some constant
$H>0$. If there is an $H\!$-cone $K_V(\gamma)\subset\R^{n+1}$ such
that $\Gamma\subset K_V(\gamma)$ and $L\subset D_{\gamma}$ if $V$ is finite,
then there is a unique $\mathcal{C}^0$ $H\!$-graph with boundary $\Gamma$ which is
$\mathcal{C}^{\infty}$ in the interior.
\end{theorem}

\newpage

\begin{center}
\includegraphics%[width=1cm,height=2cm] %[2cm,3cm][5cm,10cm]
{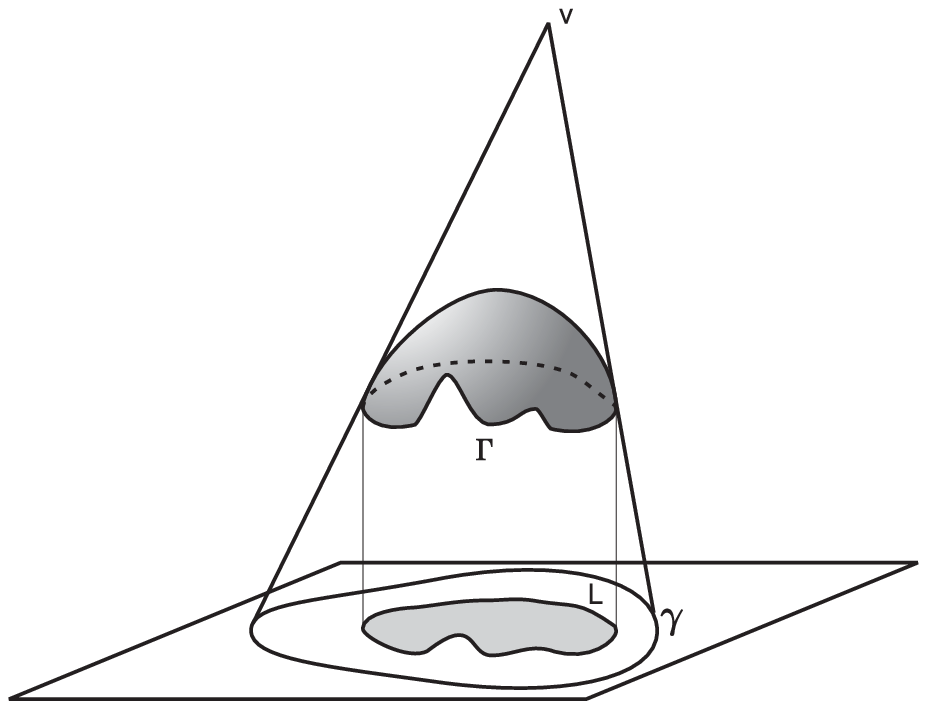}
\end{center}

\bigskip
\bigskip
\bigskip

From the above, we conclude the existence and uniqueness of $H\!$-graphs for a large class of
prescribed boundary data over domains of $\R^n$ that are not
necessarily mean convex. In fact, let $\gamma\subset\R^n$ be a
smooth compact embedded hypersurface satisfying $H_{\gamma}>H>0$. Then,  let $L$
be a smooth compact embedded hypersurface contained in the interior of the
bounded connected component of $\R^n\backslash\gamma$ with mean
curvature \mbox{$H_L\geq-H$.} Now, fix a point $P$ in the bounded connected component of
$\R^n\backslash L$ and set $\Gamma_t=K_{P+te}(\gamma)\cap C(L)$.
Clearly, there is $t_0>0$ such that $K_{P+te}(\gamma)$ is an $H\!$-cone for
any $t\geq t_0$. Therefore, there exists a unique graph over the domain enclosed by $L$
with constant mean curvature $H$ and boundary~$\Gamma_t$.
\medskip

Since a $H\!$-graph of CMC is given by a solution of a quasi-linear elliptic
second order PDE, the existence problem for $H\!$-graphs 
is equivalent to the solvability of a Dirichlet problem for the CMC equation.
Consequently, the existence of solutions is usually proved within the theory
of elliptic PDE, as is the case of the results stated above.

\bigskip

\noindent\textit{Proof of Theorem \ref{hcons smooth}}. We first consider the
case of $V$ finite. The submanifold $\Gamma$ divides $K_V(\gamma)$ into two
connected components. We smooth out the vertex of $K_V(\gamma)$ if $V$
belongs to the simply connected component that we call $G$. Then $G$ is a
graph over the domain $\Omega$ enclosed by $L$ of a function $\psi\in
C^{2,\alpha}(\bar{\Omega})$. The mean curvature function $H_{\psi}(x)$ of $G$
at $(x,\psi(x))$ belongs to $C^{1,\alpha}(\bar{\Omega})$ and satisfies
$H_{\psi}(x)\geq H$.

\newpage

Let $H_t\in C^{1,\alpha}(\bar{\Omega})$, $t\in[0,1]$, be the family of
functions
$$
H_t=(1-t)H_{\psi}+tH.
$$
Thus $H_{\psi}=H_0\ge H_t\ge H_1=H>0$. Consider the family of Dirichlet
problems
\begin{equation}\label{dir}
\left\{\begin{array}
[c]{l}
Q_t[v]=\operatorname{div}\left( \frac{\nabla v} {\sqrt{1+|\nabla v|^2}
}\right) +H_t=0\text{ in }\Omega,\text{ }v\in C^{2,\alpha}(\bar{\Omega})\\
v|_{\partial\Omega}=\varphi
\end{array}
\right.
\end{equation}
where $\varphi=\psi|_L$. The set $S$ of $t\in[0,1]$ for which (\ref{dir}) 
has a solution is non-empty since $0\in S$. Moreover, it is open by the implicit function theorem.

To prove that $S$ is closed it suffices to show the existence of a uniform
$C^1$ bound of any solution of (\ref{dir}). Since $\psi$ is a supersolution
for $Q_t$ with $\varphi=\psi|_L$ and $w=0$ is a subsolution, we have from 
the maximum principle and for any $t$ that
$$
0\leq v\leq\sup\psi\leq\langle V,e\rangle
$$
where $v\in C^{2,\alpha}(\bar{\Omega})$ is a solution of (\ref{dir}).
Therefore,
$$
\left\vert v\right\vert_0\leq\langle V,e\rangle.
$$

To estimate the gradient of $v$ at the boundary it is enough to construct a
local barrier from below in a neighborhood of $L$ in $\bar{\Omega}$ since
$\psi$ is a global barrier from above with bounded gradient (see p.\ 333 of
\cite{GT}). Then, a global barrier for the gradient follows from either Section $5$
of \cite{DHL} or Lemma $6$ of \cite{DL}.

Let $\eta$ be the unit normal vector of $L$ pointing to $\Omega$. There is
$\epsilon>0$ such that the normal exponential map
$$
E(s,y)=y+s\eta(y) \;\;\;\mbox{for}\; s\in[0,\epsilon]\;\mbox{and}\;y\in L,
$$
is a diffeomorphism from $[0,\epsilon]\times L$ onto a closed neighborhood
$\Lambda=E([0,\epsilon]\times L)$ of $L$ on $\bar{\Omega}$. We take local
coordinates $\{x_1 ,\ldots,x_n\}$ on $\Lambda$ where $x_1=s=d(x)$ is the
distance function to $L$ and the remaining are local coordinates for $L$. We
denote the corresponding coordinate vector fields by $\partial_1
=\partial_s ,\partial_2,\dots,\partial_n$. Then, the metric $\sigma
_{ij}dx_idx_{j}$ satisfies $\sigma_{11}=1$ and $\sigma_{1j}=0$ if $j\geq2$.

Choose $t\in[0,1]$. Then $Q_t[v]$ takes the form
$$
Q_t[v]=\frac{1}{A^{1/2}}\left(\sigma^{ij}-\frac{v^iv^j}{A}\right)
v_{i;j}+H_t=0
$$
where we denote $A=1+|\nabla v|^2$.

Let $w\in C^{2,\alpha}(\Lambda)$ be given by
$$
w(x)=\xi(s)+\phi(x)
$$
where $\xi(s)$ is a real function to be chosen later and $\phi(x)$ is defined
by extending $\varphi$ to $\Lambda$ as constant along the straight lines normal  to
$L$. Then, we have
$$
Q_t[w]=\frac{1}{A^{1/2}}\left(\sigma^{ij}
-\frac{(\xi^i+\phi^i)(\xi^j
+\phi^j)}{A}\right)(\xi_{i;j}+\phi_{i;j})+H_t
$$
where
$$
A=1+|\nabla\phi(x)|^2+\xi_s^2(d(x))=B+\xi_s^2.
$$

 From $|\nabla d|=1$, we have $d^id_{i;j}=0$. Since $2\langle\nabla_{\partial d}\nabla d,
\partial d\rangle=\partial d|\nabla d|^2=0$, we obtain
$$
\Delta d|_{d=s}=\sigma^{ij}\langle\nabla_{\partial_i}\nabla d, 
\partial_{j}\rangle=-\sigma^{ij}b_{i,j}(s)=-H^s
$$
where $H^s$ denotes the mean curvature of $E(L\times\{s\})$. Then,
\begin{equation}
\label{dif}Q_t[w]\geq\frac{1}{A^{1/2}}\left(\sigma^{ij} 
-\frac{(\xi^i+\phi^i)(\xi^j+\phi^j)}{A}\right) \xi_{i;j} -\Lambda|\phi|_2+H_t
\end{equation}
where $\Lambda=1/A^{1/2}$ is the largest eigenvalue in $\Omega$ of $Q_t$.

We conclude from (\ref{dif}) using
$$
\xi^i\xi^j\xi_{i;j}=\xi_s^2d^id^j(\xi_{ss}d_id_{j}+\xi
_sd_{i;j})=\xi_s^2\xi_{ss},
$$
$$
\xi^i\phi^j\xi_{i;j}=\xi_sd^i\phi^j(\xi_{ss}d_id_{j}+\xi
_sd_{i;j})=\xi_s\xi_{ss}\langle\nabla d,\nabla\phi\rangle=0,
$$
$$
\phi^i\phi^j\xi_{i;j}=\phi^i\phi^j(\xi_{ss}d_id_{j}+\xi_s
d_{i;j})=-\phi^i\phi^j\xi_sb_{i,j}(s)
$$
and
$$
\sigma^{ij}\xi_{i;j}=\Delta\xi=\xi_{ss}+\xi_s\Delta d=\xi_{ss}-\xi_sH^s
$$
that
$$
A^{3/2}Q_t[w]\geq B\xi_{ss}-AH^s\xi_s+\phi^i\phi^jb_{i,j}(s)\xi_s
-A^{3/2}\Lambda|\phi|_2+A^{3/2}H_t.
$$

We take $\xi\in C^{\infty}([0,\epsilon])$ of the form
$$
\xi(s)=\delta\ln(1+\beta s)
$$
for constants $\delta<0$ and $\beta>0$ to be determined. Thus
$$
\xi_s=\frac{\delta\beta}{1+ \beta s}<0\;\;\;\mbox{and}\;\;\;\xi_{ss} =-\frac{1}{\delta}\xi_s^2
$$
We obtain,
$$
A^{3/2}Q_t[w]\geq-B\delta^{-1}\xi_s^2-(B+\xi_s^2)H^s\xi_s
+\phi^i\phi^jb_{i,j}(s)\xi_s+(B+\xi_s^2)^{3/2}H_t-(B+\xi_s^2)|\phi|_2.
$$
Since $B=1+|\nabla\phi(x)|^2\geq1$ and $\xi_s<0$, we have
$$
A^{3/2}Q_t[w]\geq(H^s+H_t)(-\xi_s)^{3}-(B\delta^{-1}+|\phi|_2
)\xi_s^2+C\xi_s+D
$$
where the functions $C$ and $D$ depend only on the the metric of $\Lambda$ and
on the function $\varphi$ and its derivatives. Since $H^0=H_L$ and
$H_t\geq H\geq-H_L$, it follows that
$$
H^0+H_t=H_L+H_t\geq H_L+H\ge0.
$$
Therefore, at points of $L$ we have that
$$
\lim_{s\to 0}A^{3/2}Q_{t}(w)\geq-(B\delta^{-1}+|\phi|_2)(\delta\beta)^2+C\delta\beta+D.
$$
We choose $\delta$ such that $B/\delta+|\phi|_2<0$. 
Then, there is $\beta$ independent of $t$ and large enough such that
$A^{3/2}Q_t[w]\geq1$ in the neighborhood 
$\Lambda^{\prime} =E(L\times[0,\epsilon_1])$ of $L$ on $\bar{\Omega}$ for some 
$\epsilon_1\in(0,\epsilon]$.

To assure that $w$ is a local barrier from below for $Q_t$ in 
$\Lambda^{\prime}$ we have to guarantee that
\begin{equation}\label{tri}
w|_{\partial\Lambda^{\prime}}\leq v|_{\partial\Lambda^\prime}.
\end{equation}
At $L$ we have that $w=\varphi$ and (\ref{tri}) is trivially satisfied. At the
other component of the boundary of $\Lambda^{\prime}$ condition (\ref{tri}) is
verified if
$$
w(x,\epsilon_1)=\delta\ln(1+\epsilon_1\beta)+\varphi(x)\leq\delta\ln
(1+\epsilon_1\beta)+\langle V,e\rangle\leq0,
$$
and this condition is clearly satisfied by choosing $\beta$ large enough.

We have proved that $w$ is a barrier from below of $Q_t$ in $\Lambda^\prime$ 
by appropriate choices of $\delta$  and $\beta$ independent of $t$. It follows that
the $C^1$-norm of $w$ in $\Omega$ can be estimate by a bound which does not
depend on $t$. Since $w\leq v\leq\psi$ in $\Lambda$ and
$$
v|_L=\varphi=w|_L=\psi|_L
$$
by the gradient maximum principle (cf.\ Theorem 15.1 in \cite{GT}), we have
$$
\max_\Omega|\nabla v| =\max_L|\nabla v|\leq\max\{\max_L|\nabla w|
,\max_{\Omega}|\nabla\psi|\}\leq M
$$
where $M>0$ does not depend on $t$. Finally, we take $C=\max\left\{ \langle
V,e\rangle,M\right\} $. This proves the existence of a priori $C^1$
estimates for the solutions of (\ref{dir}) which implies that $S$ is closed and
concludes the proof of the theorem in the case of finite $V$. \medskip

We now use the previous case to prove the theorem for $V=\infty$, that is, to
obtain Serrin's result. Given $n\in\mathbb{N}$ and setting
$$
H_n=\frac{n}{n+1}H,
$$
there is $V_n$ high enough such that the cone
$$
J_n=\{tV_n+(1-t)p:p\in\Gamma\text{ and }t\in[0,1]\},
$$
is an $H_n$-cone. Moreover, since $H_L\geq H\geq H_n$, we may use the
previous case to assert the existence of a solution $u_n\in C^{2,\alpha
}(\bar\Omega)$ of $Q_{H_n}=0$ in $\Omega$ such that $u_n|_{\partial\Omega
}=\varphi$. We have the well-known height estimates
$$
|u_n|_0\leq\frac{2}{H_n}\le\frac{4}{H}.
$$
Moreover, a similar local barrier used above to estimate the gradient of $u_n$
from below can be used to estimate the gradient of $u_n$ from above and, as
before, this provides uniform (that is, no depending on $n)$ $C^1$ estimates
of the sequence $u_n$. From linear elliptic PDE theory this guarantees
$C^2$ compactness of $\{u_n\}$. Therefore, a subsequence of $u_n$
converges to a solution $u\in C^2$ of $Q_H=0$ such that $u|_{\partial
\Omega}=\varphi$. Finally, PDE regularity implies that $u\in C^{2,\alpha}
(\bar{\Omega})$.\qed

\bigskip

\noindent\textit{Proof of Theorem }\textit{\ref{hcons c0}}. Due to Serrin's
result for the $C^0$ case (Theorem 16.11 in \cite{GT}) it suffices to
consider the case of $V$ finite. Without loss of generality, we assume that
the domain enclosed by $\gamma$ contains the origin of $\R^n$.
Clearly, the domain $D_{\gamma_k}$ enclosed by $\gamma_k=(1+1/k)\gamma$
contains $\bar\Omega$ for any $k$. Set $V_k=(1+1/k)V$. Then, the cone
$$
K_k=K_{V_k}(\gamma_k)=(1+1/k)K_V(\gamma)
$$
is an $H_k\!$-cone for $H_k=(1+1/k)H$. Since $H_k>H$, we thus have
\begin{equation}
\label{pinva}H_L\geq-H>-H_k.
\end{equation}

Let $\psi_k\in C^0(\bar{\Omega})$ be the function whose graph is
$$
S_k:=\bar{\Omega}\times[0,+\infty)\cap K_k
$$
and $\varphi_k\in C^{2,\alpha}(L)$ the function which graph is $\partial
S_k$. We apply Perron's method on $\bar\Omega$ to the constant mean
curvature hypersurface elliptic PDE
\begin{equation}
\label{edh}Q_{H_k}[v]=\text{div}\left(\frac{\nabla v}{\sqrt{1+|\nabla v|^2}}\right) +H_k=0
\end{equation}
in $\Omega$ with boundary condition $\psi_k|_L=\varphi_k$.

It is clear that $\psi_k$ is a supersolution for (\ref{edh}). On the other
hand, we have from (\ref{pinva}) that the mean curvature of the cylinder
$L\times\R$ is strictly smaller than $-H_k$. Then, there is
$z_0\gg0$ such that the cone
$$
J_k=\{t(V-z_0e)+(1-t)p:p\in\partial S_k\text{ and }t\in[0,1]\}
$$
has mean curvature smaller than $-H_k$. If $J_k$ is the graph of the
function $\chi_k\in C^0(\bar{\Omega})$ then $\chi_k$ is a subsolution
for $Q_{H_k}$ in $\bar\Omega$ such that $\chi_k|_L =\varphi_k$. It
follows from Perron's method that
$$
v_k(x)=\sup\left\{\sigma(x):\sigma\in C^0(\bar{\Omega})\text{ is a
subsolution for }Q_{H_k}\text{ in }\Omega,\text{ }\sigma|_L =\varphi
_k\right\} ,\text{ }x\in\bar{\Omega}
$$
is in $C^{2,\alpha}(\Omega)\cap C^0(\bar{\Omega})$ and is a solution of
$Q_{H_k}=0$ in $\Omega$ such that $v_k|_L=\varphi_k$.

We have that $v_k\in C^{\infty}(\Omega)$ by interior regularity. Moreover,
since $\chi_k\leq u\leq\psi_k$ for all $k$ and $\varphi_k=\chi_k
|_{\partial\Omega}=\psi_k|_{\partial\Omega}$ converges to $\varphi$ as
$k\rightarrow\infty$, it follows that $u$ extends continuously to $\bar
{\Omega}$ and $u|_{\partial\Omega}=\varphi$.\qed

\begin{remark}\po
\emph{It is clear that we may replace in Theorems
1 and 2 the $H\!$-cone by any supersolution of $Q_H$. Observe that
$H\!$-cone like supersolutions can be constructed as follows: Consider a
star-shaped compact embedded hypersurface $\gamma\subset\R^n$ with
respect to a point $O\in D$ and take a point $V\in\R_{+}^n
\cup\{\infty\} $ in the half straight line $r$ orthogonal to $\R^n$
with origin $O$. Then, in each hyperplane $R$ containing $r$ choose a smooth
simple curve $\alpha_R$ starting in $V$ (asymptotic to $r$ if $V=\infty)$
finishing in $\gamma$ and having a single valued projection on 
$\R^n$. Moreover, the $\alpha_R$'s are chosen to depend smoothly on $R$ so
that their union as $R$ varies constitute a graph G in $\R^{n+1}$
which is smooth in $G\backslash\{V\}$ (in $G$ if $V=\infty)$. Given $H\geq0$,
it is clear that we may construct such a $G$ to be an \mbox{$H-$supersolution}
for $Q_H$ on a domain containing $O$. }
\end{remark}

{\begin{tabbing}
\indent \= Marcos Dajczer\hspace{25,5ex} Jaime Ripoll\\
\> IMPA \hspace{35ex}Instituto de Matematica\\
\> Estrada Dona Castorina, 110\hspace{12ex}
Univ. Federal do Rio Grande do Sul\\
\> 22460-320 -- Rio de Janeiro -- RJ
\hspace{8ex} Av. Bento Gon\c calves 9500\\
\> Brazil\hspace{35ex} 91501-970 -- Porto Alegre -- RS\\
\> marcos@impa.br\hspace{24,6ex} Brazil\\
\>\hspace{41ex} jaime.ripoll@ufrgs.br
\end{tabbing}}

\end{document}